\documentclass{amsart}
\usepackage{graphicx}
\usepackage{amssymb}
\usepackage{amsfonts}
\usepackage{bm}
\swapnumbers
\newtheorem{theorem}{Theorem}[section]

\theoremstyle{definition}

\newtheorem{assumption}[theorem]{Assumption}

\numberwithin{equation}{section}
 \theoremstyle{plain}
 
 \numberwithin{equation}{section} 
 \numberwithin{figure}{section} 
 \theoremstyle{plain}
 \theoremstyle{remark}
 \newtheorem*{acknowledgement*}{Acknowledgement}

\newcommand{\cA}{{\mathcal A}}

\newcommand{\cE}{{\mathcal E}}
\newcommand{\cF}{{\mathcal F}}
\newcommand{\cG}{{\mathcal G}}
\newcommand{\cH}{{\mathcal H}}

\newcommand{\cN}{{\mathcal N}}

\newcommand{\Om}{{\Omega}}
\newcommand{\om}{{\omega}}
\newcommand{\ve}{{\varepsilon}}
\newcommand{\del}{{\delta}}
\newcommand{\Del}{{\Delta}}
\newcommand{\gam}{{\gamma}}
\newcommand{\Gam}{{\Gamma}}

\newcommand{\sig}{{\sigma}}
\newcommand{\al}{{\alpha}}
\newcommand{\be}{{\beta}}

\newcommand{\la}{{\lambda}}
\newcommand{\La}{{\Lambda}}

\newcommand{\bbN}{{\mathbb N}}

\newcommand{\bbR}{{\mathbb R}}

\newcommand{\bbZ}{{\mathbb Z}}
\newcommand{\bbI}{{\mathbb I}}

\begin{document}
\title[]{The strong Borel-Cantelli property\\
in conventional and nonconventional setups}%
 \vskip 0.1cm
 \author{ Yuri Kifer\\
\vskip 0.1cm
 Institute  of Mathematics\\
Hebrew University\\
Jerusalem, Israel}%
\address{
Institute of Mathematics, The Hebrew University, Jerusalem 91904, Israel}
\email{ kifer@math.huji.ac.il}%

\thanks{ }
\dedicatory{ }
 \date{\today}
\begin{abstract}\noindent
We study the strong Borel-Cantelli property both for events and for shifts on sequence
spaces considering both a conventional and a nonconventional setups. Namely, under certain conditions on events
$\Gam_1,\Gam_2,...$ we show that with probability one
\[
(\sum_{n=1}^N\prod_{i=1}^\ell P(\Gam_{q_i(n)}))^{-1}\sum_{n=1}^N\prod_{i=1}^\ell\bbI_{\Gam_{q_i(n)}}\to 1\,\,\mbox{as}\,\, N\to\infty
\]
where $q_i(n),\, i=1,...,\ell$ are integer valued functions satisfying certain assumptions and $\bbI_\Gam$ denotes the indicator
of $\Gam$. When $\ell=1$ (called the conventional setup) this convergence can be established under $\phi$-mixing conditions while
when $\ell>1$ (called a nonconventional setup) the stronger $\psi$-mixing condition is required. These results are extended to
shifts $T$ of sequence spaces where $\Gam_{q_i(n)}$ is replaced by $T^{-q_i(n)}C^{(i)}_n$ where $C_n^{(i)},i=1,...,\ell,\, n\geq 1$ is a sequence of cylinder sets.
As an application we study the asymptotical behavior of maximums of certain logarithmic  distance functions and of ( multiple) hitting times of shrinking cylinders.
\end{abstract}
\maketitle
\markboth{Yu.Kifer}{Borel-Cantelli property}
\renewcommand{\theequation}{\arabic{section}.\arabic{equation}}
\pagenumbering{arabic}

\section{Introduction}\label{sec1}\setcounter{equation}{0}

The classical second Borel--Cantelli lemma states that if $\Gam_1,\Gam_2,...$ is a sequence of independent
events such that
\begin{equation}\label{1.1}
\sum_{n=1}^\infty P(\Gam_n)=\infty
\end{equation}
then with probability one infinitely many of events $\Gam_i$ occur, i.e.
\begin{equation}\label{1.2}
\sum_{n=1}^\infty\bbI_{\Gam_n}=\infty\quad\mbox{almost surely (a.s.)}
\end{equation}
where $\bbI_\Gam$ is the indicator of a set (event) $\Gam$.

There is a long list of papers, starting probably with \cite{Le}, providing conditions which replace the
independency by a weaker assumption and still yield (\ref{1.2}) (see, for instance, \cite{DMR} and references
there). On the other hand, it was shown in Theorem 3 of \cite{Ph} that under $\phi$-mixing with a summable coefficient
$\phi$ the condition (\ref{1.1}) yields the stronger version of the second Borel--Cantelli lemma in the form
\begin{equation}\label{1.3}
\frac {S_N}{\cE_N}\to 1\,\,\,\mbox{almost surely (a.s.) as}\,\, N\to\infty
\end{equation}
where $S_N=\sum_{n=1}^N\bbI_{\Gam_n}$ and $\cE_N=\sum_{n=1}^NP(\Gam_n)$.

The same paper \cite{Ph} started another line of research, known now under the name dynamical Borel--Cantelli lemmas,
where (\ref{1.3}) is proved for $S_N=\sum_{n=1}^N\bbI_{\Gam_n}\circ T^n$ where $T$ is a measure preserving transformation
on a probability space $(\Om,P)$ and $\Gam_n,\, n\geq 1$ is a sequence of measurable sets. For such $S_N$'s the convergence
(\ref{1.3}) was proved, in particular, for the Gauss map $Tx=\frac 1x$ (mod 1), $x\in(0,1]$ preserving the Gauss measure
$P(\Gam)=\frac 1{\ln 2}\int_\gam\frac {dx}{1+x}$. This line of research became quite popular in the last two decades.
In particular, \cite{CK} proves (\ref{1.3}) in the dynamical setup considering $T$ being the so called subshift of finite
type on a sequence space  where $\Gam_n,\, n\geq 1$ is a sequence of cylinders while another series of papers dealt with
uniformly and non-uniformly hyperbolic dynamical systems as a transformation $T$ and with geometric balls as $\Gam_n$'s
(see, for instance, \cite{Gu}, \cite{HNPV} and references there).

In this paper we consider, in particular, "nonconventional" extensions of some of the above results aiming to prove that under
 certain conditions (\ref{1.3}) holds true with $S_N=\sum_{n=1}^N(\prod_{i=1}^\ell\bbI_{\Gam_{q_i(n)}})$ and $\cE_N=\sum_{n=1}^NP(\Gam_{q_i(n)})$
 where $q_i(n),\, i=1,...,\ell$ functions taking on positive integer values on positive integers and satisfying certain assumptions valid, in particular,
 for polynomials with integer coefficients. When $\ell=1$ (conventional setup) the $\phi$-mixing with a summable coefficient $\phi$
 suffices for our result, while for $\ell>1$ we have to impose stronger $\psi$-mixing conditions.

 In the dynamical systems setup we consider $S_N=\sum_{n=1}^N(\prod_{i=1}^\ell\bbI_{C_n^{(i)}}\circ T^{q_i(n)})$ and
 $\cE_N=\sum_{n=1}^N\prod_{i=1}^\ell P(C_n^{(i)})$ where $T$ is the left shift on a sequence space $\cA^{\bbZ}$ with a finite or countable
 alphabet while $C_n^{(i)},\, i=1,...,\ell,\, n\geq 1$ is a sequence of cylinder sets. As an application we study the asymptotic behaviors of expressions
 $M_N=\max_{1\leq n\leq N}(\min_{1\leq i\leq\ell}\Phi_{\tilde\om^{(i)}}\circ T^{q_i(n)})$ where $\Phi_{\tilde\om}(\om)=-\ln(d(\om,\tilde\om))$,
 $\om,\tilde\om\in\cA^\bbN$ and $d(\cdot,\cdot)$ is the natural distance on the sequence space.

 Our results extend some of the previous work in the
 following aspects. First, the strong Borel--Cantelli property in the nonconventional setup $\ell>1$ was not studied
 before at all. Secondly, even in the conventional setup $\ell=1$ considering rather general functions $q(n)=q_1(n)$
 in place of just $q(n)=n$ seems to be new, as well. Thirdly, we extend for shifts some of the results from \cite{CK}
 considering sequence spaces with countable alphabets and $\phi$-mixing invariant measures rather than just subshifts
 of finite type with Gibbs measures which are exponentially fast $\psi$-mixing (see \cite{Bo}). This allows to apply our
 results, for instance, to Gibbs-Markov maps and to Markov chains with a countable state space satisfying the Doeblin condition
 since both examples are exponentially fast $\phi$-mixing, see \cite{MN} and \cite{Br}, respectively.

 In the next section we will formulate precisely our setups and assumptions and state our main results. In Section \ref{sec3} we
 will prove the strong Borel--Cantelli property for events under the $\phi$-mixing condition in the conventional setup
 $\ell=1$ and under $\psi$-mixing condition in the nonconventional setup $\ell>1$. In Sections \ref{sec4} and \ref{sec5} we extend
 the strong Borel--Cantelli property to shifts under the $\phi$-mixing when $\ell=1$ and under $\psi$-mixing when $\ell>1$, respectively.
 In Section \ref{sec6} we exhibit applications to the asymptotic behaviors of maximums along shifts of logarithmic distance functions while in the last Section \ref{sec7} we apply the strong
 Borel--Cantelli property to derive the asymptotics of multiple hitting times of shrinking cylinder
  sets.

 \section{Preliminaries and main results}\label{sec2}\setcounter{equation}{0}

We start with a probability space $(\Om,\cF,P)$ and a two parameter family of
$\sig$-algebras $\cF_{mn}$ indexed by pairs of integers $-\infty\leq m\leq n\leq\infty$ and
such that $\cF_{mn}\subset\cF_{m'n'}\subset\cF$ if $m'\leq n\leq n'$. Recall that the
$\phi$ and $\psi$ dependence coefficient between two $\sig$-algebras $\cG$ and $\cH$ can be
written in the form
(see \cite{Br}),
\begin{eqnarray}\label{2.1}
&\phi(\cG,\cH)=\sup_{\Gam\in\cG,\,\Del\in\cH}\{|\frac {P(\Gam\cap\Del)}{P(\Gam)}-P(\Del)|,\, P(\Gam)\ne 0\}\\
&=\frac 12\sup\{\| E(g|\cG)-Eg\|_{L^\infty}:\, g\,\,\mbox{is $\cH$-measurable and}\,\,\| g\|_{L^\infty}\leq 1\}
\nonumber\end{eqnarray}
and
\begin{eqnarray}\label{2.2}
&\psi(\cG,\cH)=\sup_{\Gam\in\cG,\,\Del\in\cH}\{|\frac {P(\Gam\cap\Del)}{P(\Gam)P(\Del)}-1|,\, P(\Gam)P(\Del)\ne 0\}\\
&=\frac 12\sup\{\| E(g|\cG)-Eg\|_{L^\infty}:\, g\,\,\mbox{is $\cH$-measurable and}\,\, E|g|\leq 1\},
\nonumber\end{eqnarray}
respectively. The $\phi$-dependence (mixing) and the $\psi$-dependence (mixing) in the family $\cF_{mn}$ is measured
by the coefficients
\begin{equation}\label{2.3}
\phi(k)=\sup_m\phi(\cF_{-\infty,m},\cF_{m+k,\infty})\,\,\mbox{and}\,\,\psi(k)=\sup_m\psi(\cF_{-\infty,m},\cF_{m+k,\infty}),
\end{equation}
respectively, where $k=0,1,2,...$. The probability measure $P$ is called $\phi$-mixing or $\psi$-mixing with respect to the family
of $\sig$-algebras $\cF_{mn}$ if $\phi(n)\to 0$ or $\psi(1)<\infty$ and $\psi(n)\to 0$ as $n\to\infty$, respectively.

Our setup includes also functions $q_1(n),\, q_2(n),...,q_\ell(n)$ with $\ell\geq 1$ taking on nonnegative integer values on
integers $n\geq 0$ and satisfying
\begin{assumption}\label{ass2.1} There exists a constant $K>0$ such that

(i) for any $i\ne j$, $1\leq i,j\leq\ell$ and every integer $k$ the number of integers $n\geq 0$ satisfying at least one of
the equations
\begin{equation}\label{2.4}
q_i(n)-q_j(n)=k\quad\mbox{and}\quad q_i(n)=k
\end{equation}
does not exceed $K$ (when $\ell=1$ only the second equation in (\ref{2.4}) should be taken into account);

(ii) the cardinality of the set $\cN$ of all pairs $n>m\geq 0$ satisfying
\begin{equation}\label{2.5}
\max_{1\leq i\leq\ell}q_i(n)\leq\max_{1\leq i\leq\ell}q_i(m)
\end{equation}
does not exceed $K$.
\end{assumption}

Observe that Assumption \ref{ass2.1} is satisfied if $q_i,\, i=1,...,\ell$ are essentially distinct nonconstant
polynomials (i.e. $|q_i(n)-q_j(n)|\to\infty$ as $n\to\infty$ for any $i\ne j$) with integer coefficients taking on
nonnegative values on nonnegative integers. Indeed, $q_i(n)-q_j(n)$ and $q_i(n)$ are nonconstant polynomials,
and so the number of $n$'s solving one of equations in (\ref{2.4}) is bounded by the degree of the corresponding polynomial.
In order to show that (\ref{2.5}) can hold true in the polynomial case only for finitely many pairs $m<n$ observe that
there exists $n_0\geq 1$ such that all polynomials $q_1(n),\, q_2(n),...,q_\ell(n)$ are strictly increasing on $[n_0,\infty)$.
Hence, if $n>m\geq n_0$ then (\ref{2.5}) cannot hold true. If $0\leq m<n_0$ and $n\geq n_0$ then there exists $n_1\geq n_0$ such
that for all $n\geq n_1$ (\ref{2.5}) cannot hold true, as well. The remaining case $0\leq m<n_0$ and $0\leq n<n_1$ concerns less
than $n_0n_1$ pairs $m<n$.

Next, we will state our result concerning sequences of events. Let $\Gam_1,\Gam_2,...\in\cF$ be a sequence of events and each
$\sig$-algebra $\cF_{mn},\, 1\leq m\leq n<\infty$ be generated by the events $\Gam_m,\Gam_{m+1},...,\Gam_n$. Set also
$\cF_{mn}=\cF_{1n}$ for $-\infty\leq m\leq 0$ and $n\geq 1$, $\cF_{mn}=\{\emptyset,\Om\}$ for $m,n\leq 0$ and $\cF_{m,\infty}=
\sig\{\Gam_m,\Gam_{m+1},...\}$. Set
\begin{equation}\label{2.6}
S_N=\sum_{n=1}^N(\prod_{i=1}^\ell\bbI_{\Gam_{q_i(n)}})\quad\mbox{and}\quad\cE_N=\sum_{n=1}^N(\prod_{i=1}^\ell P(\Gam_{q_i(n)})).
\end{equation}
\begin{theorem}\label{thm2.2}
Let $\phi$ and $\psi$ be dependence coefficients defined by (\ref{2.3}) for the above $\sig$-algebras $\cF_{mn}$. Assume that
$\phi(n),\, n\geq 0$ is summable in the case $\ell=1$ and $\psi(n),\, n\geq 0$ is summable in the case $\ell>1$. Suppose that
the functions $q_1(n),...,q_\ell(n)$ satisfy Assumption \ref{ass2.1}(i) and
\begin{equation}\label{2.7}
\cE_N\to\infty\quad\mbox{as}\quad N\to\infty.
\end{equation}
Then, with probability one,
\begin{equation}\label{2.8}
\lim_{N\to\infty}\frac {S_N}{\cE_N}=1\quad\mbox{as}\quad N\to\infty.
\end{equation}
\end{theorem}

Next, we will present our results concerning shifts. Here $\Om=\cA^\bbZ$ is the space of sequences $\om=(...,\om_{-1},\om_0,\om_1,...)$
with terms $\om_i$ from a finite or countable alphabet $\cA$ which is not a singleton with the index $i$ running along integers (or along
natural numbers $\bbN$ which can also be considered requiring very minor modifications). We assume that the basic $\sig$-algebra $\cF$ is
generated by all cylinder sets while the $\sig$-algebras $\cF_{mn},\, n\geq m$ are generated by the cylinder sets of the form
$\{\om=(\om_i)_{-\infty<i<\infty}:\,\om_i=a_i\,\,\mbox{ for}\,\,\, m\leq i\leq n\}$ for some $a_m,a_{m+1},...,a_n\in\cA$. The setup includes also the left shift $T:\Om\to\Om$ acting by $(T\om)_i=\om_{i+1}$ and a $T$-invariant probability measure $P$ on $(\Om,\cF)$, i.e. $P(T^{-1}\Gam)=P(\Gam)$
for any measurable $\Gam\subset\Om$. In this setup $\phi$ and $\psi$-dependence coefficients defined by (\ref{2.3}) will be considered with
respect to the family of $\sig$-algebras $\cF_{mn},\, m\leq n$ defined above. Without loss of generality we assume that the probability of each
1-cylinder $[a]=\{\om=(\om_i)_{i\in\bbZ}:\, \om_0=a$ is positive, i.e. $P([a])>0$ for any $a\in\cA$, and since $\cA$ is not a singleton we have
also that $\sup_{a\in\cA}P([a])<1$.

Each cylinder $C$ is defined on an interval of integers $\La=[l,r]$, $l\leq r$, i.e. $C=\{\om=(\om_i)_{-\infty<i<\infty}:\,\om_i=a_i,\, i=l,l+1,...,r\}$
for some $a_l,...,a_r\in\cA$. Given a constant $D>0$ call an interval of integers $\La_1=[l_1,r_1]$ to be right $D$-nested in the interval of
integers $\La_2=[l_2,r_2]$ if $[l_1,r_1]\subset (-\infty,r_2+D)$, i.e. $r_1<r_2+D$. Such an interval $\La_1$ will be called $D$-nested in $\La_2$
if $[l_1,r_1]\subset (l_2-D,r_2+D)$. The latter notion was used also in \cite{CK}.

Let $C_n^{(j)},\, j=1,...,\ell,\, n=1,2,...$ be a sequence of cylinder sets defined on intervals of integers $\La_n,\, n=1,2,...$ so that $C_n^{(j)},\, j=1,...,\ell$ are defined on $\La_n$ for each $n\geq 1$. Set
\begin{equation}\label{2.9}
S_N=\sum_{n=1}^N(\prod_{i=1}^\ell\bbI_{C_n^{(i)}}\circ T^{q_i(n)})\quad\mbox{and}\quad\cE_N=\sum_{n=1}^N
\prod_{i=1}^\ell P(C_n^{(i)}).
\end{equation}
\begin{theorem}\label{thm2.3} Suppose that the functions $q_1(n),...,q_\ell(n)$ satisfy Assumption \ref{ass2.1} and
\begin{equation}\label{2.10}
\cE_N\to\infty\quad\mbox{as}\quad N\to\infty.
\end{equation}
Let $C_n^{(j)},\, j=1,...,\ell,\, n\geq 1$ be a sequence of cylinder sets defined on intervals $\La_n\subset\bbZ$ as
described above and $D>0$ be a constant.

(i) If $\ell=1$ assume that the $\phi$-dependence coefficient is summable and that for all $m<n$ the interval
$\La_m$ is right $D$-nested in $\La_n$. Then, with probability one,
\begin{equation}\label{2.11}
\lim_{N\to\infty}\frac {S_N}{\cE_N}=1\quad\mbox{as}\quad N\to\infty.
\end{equation}
(ii) If $\ell>1$ assume that the $\psi$-dependence coefficient is summable and that for all $m<n$ the interval
$\La_m$ is $D$-nested in $\La_n$. Then with probability one (\ref{2.11}) holds true, as well.
\end{theorem}

As in most papers on the strong Borel--Cantelli property both Theorems \ref{thm2.2} and \ref{thm2.3} rely on the following basic result.
\begin{theorem}\label{thm2.4}
Let $\Gam_1,\Gam_2,...$ be a sequence of events such that for any $N\geq M\geq 1$,
\begin{equation}\label{2.12}
\sum_{m,n=M}^N(P(\Gam_m\cap\Gam_n)-P(\Gam_m)P(\Gam_n))\leq c\sum_{n=M}^NP(\Gam_n)
\end{equation}
where a constant $c>0$ does not depend on $M$ and $N$. Then for each $\ve>0$ almost surely
\begin{equation}\label{2.13}
S_N=\cE_N+O(\cE_N^{1/2}\log^{\frac 32+\ve}\cE_N)
\end{equation}
where
\[
S_N=\sum_{n=1}^N\bbI_{\Gam_n}\quad\mbox{and}\quad\cE_N=\sum_{n=1}^NP(\Gam_n).
\]
In particular, if
\[
\cE_N\to\infty\quad\mbox{as}\quad N\to\infty
\]
then with probability one
\[
\lim_{N\to\infty}\frac {S_N}{\cE_N}=1\quad\mbox{as}\quad N\to\infty.
\]
\end{theorem}

This result (as well as the part of Theorem \ref{thm2.2} for $\ell=1$ and $q_1(n)=n$) appears
already in Theorem 3 from \cite{Ph} and in a slightly more general (analytic) form it is proved as Lemma 10
in \S 7 of Ch.1 from \cite{Sp}. Both sources refer to \cite{Sch} as the origin of this result.

We observe that Theorem \ref{thm2.3} extends Theorem 2.1 from \cite{CK} in several directions. First, for
$\ell=1$ we prove the result for arbitrary $\phi$-mixing probability measures with a summable coefficient
$\phi$ on a shift space with a countable alphabet and not just for subshifts of finite type with Gibbs
measures. Secondly, the case $\ell>1$ and rather general functions $q_i(n)$ in place of just $\ell=1$
and $q_1(n)=n$ were not considered before both in the setups of Theorem \ref{thm2.2} and \ref{thm2.3}.

A direct application of Theorem \ref{thm2.3} yields corresponding strong Borel--Cantelli property for
dynamical systems which have symbolic representations by means of finite or countable partitions, for instance,
hyperbolic dynamical systems (see, for instance, \cite{Bo}) where sequences of cylinders in Theorem \ref{thm2.3}
should be replaced by corresponding sequences of elements of joins of iterates of the partition. By a slight
modification (just by considering cylinder sets defined on intervals of nonnegative integers only) Theorem \ref{thm2.3}
remains valid for one-sided shifts and then it can be applied to noninvertible dynamical systems having a symbolic
representation via their finite or countable partitions such as expanding transformations, the Gauss map of the
interval and more general transformations generated by $f$-expansions (see \cite{He}).

In Section \ref{sec6} we apply Theorem \ref{thm2.3} to some limiting problems obtaining a symbolic version
of results from \cite{HNT} which dealt with dynamical systems on $\bbR^d$ or manifolds and not with shifts. Namely, in the
setup of Theorem \ref{thm2.3} introduce the distance between $\om=(\om_i)_{i\in\bbZ}$ and $\tilde\om=(\tilde\om_i)_{i\in\bbZ}$
from $\Om$ by
\begin{equation}\label{2.14}
d(\om,\tilde\om)=\exp(-\gam\min\{ i\geq 0:\,\om_i\ne\tilde\om_i\,\,\mbox{or}\,\,\om_{-i}\ne\tilde\om_{-i}\}),\,\,\gam>0.
\end{equation}
Set
\begin{eqnarray}\label{2.15}
&\Phi_{\tilde\om}(\om)=-\ln(d(\om,\tilde\om))\,\,\mbox{for}\,\,\om,\tilde\om\in\Om\,\,\mbox{and}\\
&M_{N,\bm{\tilde\om}}(\om)=M_{N,\tilde\om^{(1)},...,\tilde\om^{(\ell)}}=
\max_{1\leq n\leq N}\min_{1\leq i\leq\ell}(\Phi_{\tilde\om^{(i)}}\circ T^{q_i(n)}(\om))\nonumber
\end{eqnarray}
for some fixed $\ell$-tuple $\bm{\tilde\om}=(\tilde\om^{(1)},...,\tilde\om^{(\ell)}),\,\tilde\om^{(i)}\in\Om,\,
i=1,...,\ell$.
\begin{theorem}\label{thm2.5} Assume that the entropy of the partition into 1-cylinders is finite, i.e.
\begin{equation}\label{2.16}
-\sum_{a\in\cA}P([a])\ln P([a])<\infty.
\end{equation}
Then, under the conditions of Theorem \ref{2.3} for almost all $\tilde\om^{(1)},...,\tilde\om^{(\ell)}\in\Om$
 with probability one,
\begin{equation}\label{2.17}
\frac {M_{N,\tilde\om^{(1)},...,\tilde\om^{(\ell)}}}{\ln N}\to\frac \gam{2\ell h}\,\,\,\mbox{as}\,\,\, N\to\infty
\end{equation}
where $h$ is the Kolmogorov--Sinai entropy of the shift $T$ on the probability space $(\Om,\cF,P)$ and, as in Theorem \ref{thm2.3},
if $\ell=1$ we assume only $\phi$-mixing with a summable coefficient $\phi$ and if $\ell>1$ we assume $\psi$-mixing with a summable
coefficient $\psi$ (and in both cases $h>0$ by Lemma 3.1 in \cite{KY} and Lemma 3.1 in \cite{KR}).
\end{theorem}

In Section \ref{sec7} we demonstrate another application of Theorem \ref{thm2.3} deriving the
asymptotical behavior of multiple hitting times of shrinking cylinders. Namely, set
\[
\tau_{C_n(\tilde\om)}=\min\{ k\geq 1:\,\prod_{i=1}^\ell\bbI_{C_n(\tilde\om)}\circ T^{q_i(k)}(\om)=1\}
\]
where $\om,\tilde\om\in\Om$ and $C_n(\om)=\{\om=(\om_i)_{i\in\bbZ}\in\Om:\,\om_i=\tilde\om_i$ provided
$|i|\leq n\}$.

\begin{theorem}\label{thm2.6} Assume that (\ref{2.16}) holds true. Then under the conditions of
Theorem \ref{thm2.3} for $P\times P$-almost all pairs $(\om,\tilde\om)\in\Om\times\Om$,
\begin{equation}\label{2.18}
\lim_{n\to\infty}\frac 1n\ln\tau_{C_n(\tilde\om)}(\om)=2\ell h.
\end{equation}
\end{theorem}

We observe that (\ref{2.18}) was proved in \cite{KR} under the $\psi$-mixing assumption assuming
additionally stronger conditions than here while the $\phi$-mixing case was not treated there at all.
The proof of Theorem \ref{thm2.6} here is different from \cite{KR} as it relies on the Borel--Cantelli
 lemma and the strong Borel--Cantelli property which is an adaptation to our symbolic (and nonconventional) setup of proofs from \cite{Ga1} and \cite{Ga2}. We note that both Theorem
  \ref{thm2.5} and Theorem \ref{thm2.6} remain valid (with essentially the same proof) for one sided
  shifts just by deleting 2 in (\ref{2.17}) and (\ref{2.18}).

\section{Proof of Theorem 2.2}\label{sec3}\setcounter{equation}{0}
\subsection{The case $\ell=1$}\label{subsec3.1}
Let $N\geq M$ and fix an $m$ between $M$ and $N$. By Assumption \ref{ass2.1} for each $k$ there
exists at most $K$ of integers $n$ such that $q(n)-q(m)=k$ where $q(n)=q_1(n)$. If $q(n)-q(m)=k\geq 1$
then by the definition of the $\phi$-dependence coefficient
\begin{equation}\label{3.1}
|P(\Gam_{q(m)}\cap\Gam_{q(n)})-P(\Gam_{q(m)})P(\Gam_{q(n)})|\leq\phi(k)P(\Gam_{q(m)}).
\end{equation}
Hence,
\begin{equation}\label{3.2}
\sum_{N\geq n\geq M,\, q(n)>q(m)}
|P(\Gam_{q(m)}\cap\Gam_{q(n)})-P(\Gam_{q(m)})P(\Gam_{q(n)})|\leq KP(\Gam_{q(m)})\sum_{k=1}^\infty\phi(k).
\end{equation}
Since the coefficient $\phi$ is summable and that similar inequalities hold true when $q(m)>q(n)$
we conclude that the condition (\ref{2.12}) of Theorem \ref{thm2.4} is satified with $\Gam_{q(n)}$ in place of
$\Gam_n,\, n=1,2,...$ there, and so (\ref{2.8}) follows in the case $\ell=1$ assuming (\ref{2.7}). \qed

\subsection{The case $\ell>1$}.
We start with the following counting arguments concerning the functions $q_i,\, i=1,...,\ell$ satisfying
Assumption \ref{ass2.1}. Introduce
\[
q(n)=\min_{1\leq i\ne j\leq\ell}|q_i(n)-q_j(n)|.
\]
By Assumption \ref{ass2.1}(i) for each pair $i\ne j$ and any $k$ there exists at most $K$ nonnegative integers
$n$ such that $q_i(n)-q_j(n)=k$, and so
\begin{equation}\label{3.3}
\#\{ n>0:\, q(n)=k\}<K\ell^2
\end{equation}
where $\#$ stands for "the number of ...". We will need also the following semi-metric between integers $k,l>0$,
\[
\del(k,l)=\min_{1\leq i,j\leq\ell}|q_i(k)-q_j(l)|.
\]
It follows from Assumption \ref{ass2.1}(i) that for any integers $m>0$ and $k\geq 0$,
\begin{equation}\label{3.4}
\#\{ n>0:\, \del(m,n)=k\}<2K^2\ell^2.
\end{equation}
Indeed, the number of $m$'s such that $q_j(m)=q_i(n)-k$ for a fixed $i,j,n$ and $k$ does not exceed $K$
by Assumption \ref{ass2.1}(i) and (\ref{3.4}) follows since $1\leq i,j\leq\ell$.

In order to prove Theorem \ref{thm2.2} for $\ell>1$ we will estimate first
\begin{equation}\label{3.5}
|E(X_mX_n)-EX_mEX_n|=|P(\cap_{i=1}^\ell(\Gam_{q_i(m)}\cap\Gam_{q_i(n)}))-P(\cap_{i=1}^\ell\Gam_{q_i(m)})P(\cap_{i=1}^\ell\Gam_{q_i(n)})|
\end{equation}
where $m,n>0$ and $X_k=\prod_{i=1}^\ell\bbI_{\Gam_{q_i(k)}}$. If $\del(m,n)=k\geq 1$ then by Lemma 3.3 in \cite{KR}
and the definition of the $\psi$-dependence coefficient
\begin{equation}\label{3.6}
|E(X_mX_n)-EX_mEX_n|\leq 2^{2\ell+2}\psi(k)(2-(1+\psi(k))^\ell)-2EX_mEX_n
\end{equation}
where we assume, in fact, that $k$ is large enough so that $\psi(k)<2^{1/\ell}-1$. Thus, let $k_0=\min\{ k:\,\psi(k)<2^{1/\ell}-1\}$.
Then by (\ref{3.4}) and (\ref{3.6}),
\begin{equation}\label{3.7}
\sum_{N\geq n\geq M}|E(X_mX_n)-EX_mEX_n|\leq cEX_m
\end{equation}
where
\[
c=2K^2\ell^2\big(1+2^{2\ell+2}(2-(1+\psi(k_0))^\ell)^{-2}\sum_{k=k_0}^\infty\psi(k)\big)
\]
where we took into account that
\[
|EX^2_m-(EX_m)^2|\leq EX_m.
\]
Summing in (\ref{3.7}) in $m$ between $M$ an $N$ we obtain the condition (\ref{2.12}) of Theorem \ref{thm2.4}
with $\cap_{i=1}^\ell\Gam_{q_i(n)}$ in place of $\Gam_n$ there. Hence if
\begin{equation}\label{3.8}
\sum_{n=1}^\infty P(\cap_{i=1}^\ell\Gam_{q_i(n)})=\infty
\end{equation}
then Theorem \ref{thm2.4} yields that with probability one
\begin{equation}\label{3.9}
\frac {S_N}{\tilde\cE_N}\to 1\quad\mbox{as}\quad N\to\infty
\end{equation}
where $\tilde\cE_N=\sum_{n=1}^NP(\cap_{i=1}^\ell\Gam_{q_i(n)})$.

Since we assume (\ref{2.10}) and not (\ref{3.8}), it remains to show that under our conditions,
\begin{equation}\label{3.10}
\frac {\tilde\cE_N}{\cE_N}\to 1\quad\mbox{as}\quad N\to\infty.
\end{equation}
By Lemma 3.2 from \cite{KR} we obtain when $q(n)=k\geq 1$ that
\begin{equation}\label{3.11}
|P(\cap_{i=1}^\ell\Gam_{q_i(n)})-\prod_{i=1}^\ell P(\Gam_{q_i(n)})|\leq ((1+\psi(k))^\ell-1)\prod_{i=1}^\ell P(\Gam_{q_i(n)}).
\end{equation}
For $q(n)=0$ we estimate the left hand side of (\ref{3.11}) just by 1. Hence, by (\ref{3.3}),
\begin{eqnarray}\label{3.12}
&|\tilde\cE_N-\cE_N|\leq K\ell^2+\sum_{n=1,q(n)\geq 1}^N\big(((1+\psi(q(n)))^\ell-1)\prod_{i=1}^\ell P(\Gam_{q_i(n)})\big)\\
&\leq K\ell^2+\sum_{n=1,q(n)\geq 1}^N((1+\psi(q(n)))^\ell-1)\nonumber\\
&\leq K\ell^2+K\ell^2\sum_{n=1}^\infty((1+\psi(q(n)))^\ell-1)\leq C<\infty
\nonumber\end{eqnarray}
for some constant $C>0$, since the coefficient $\psi$ is summable. Dividing (\ref{3.12}) by $\cE_N$ and taking into account (\ref{2.10}) we
obtain (\ref{3.10}) and complete the proof of Theorem \ref{thm2.2}.  \qed

\section{Proof of Theorem 2.3($i$)}\label{sec4}\setcounter{equation}{0}

Here $\ell=1$, and so we set $C_n=C_n^{(1)}$ and $q(n)=q_1(n)$.
Consider cylinder sets $C_m$ and $C_n$, $1\leq m<n$ defined on intervals of integers
$\La_m=[l_m,r_m]$ and $\La_n=[l_n,r_n]$ with $\La_m$ right $D$-nested in $\La_n$
implying that $r_m<r_n+D$. Let $k=q(n)-q(m)$. By Assumption \ref{ass2.1}(i)
for each $m$ and $k$ this equality can hold true only for at most $K$ of $n$'s and by Assumption \ref{ass2.1}(ii)
for no more than $K$ of $n$'s we may have $q(n)\leq q(m)$. Next, we can write
\begin{equation}\label{4.1}
r_n+q(n)>r_m+q(m)+k-D.
\end{equation}

Assume first that
\begin{equation}\label{4.2}
l_n+q(n)\leq r_m+q(m)\quad\mbox{and}\quad r_n+q(n)>r_m+q(m).
\end{equation}
Let $C_n=[a_{l_n},a_{l_n+1},...,a_{r_n}]$ and $\hat C_{m,n}=[a_{t_{m,n}},a_{t_{m,n}+1},...,a_{r_n}]$ where we assume that
$r_n>l_n$,
\begin{eqnarray*}
&t_{m,n}=s_{m,n}+[\frac 12(r_n-s_{m,n}+1)]\quad\mbox{and}\\
& a_{m,n}=l_n+(r_m+q(m)-l_n-q(n))+1=r_m+q(m)-q(n)+1.
\end{eqnarray*}
It follows that
\begin{equation}\label{4.3}
r_n-t_{m,n}+1\geq [\frac 12(k-D)]\quad\mbox{and}\quad t_{m,n}+q(n)-r_m+q(m)\geq [\frac 12(k-D)]-1.
\end{equation}
Assuming that $k\geq D+4$ we obtain by the definition of the $\phi$-dependence coefficient that
\begin{eqnarray}\label{4.4}
&P(T^{-q(m)}C_m\cap T^{-q(n)}C_n)\leq P(T^{-q(m)}C_m\cap T^{-q(n)}\hat C_{m,n})\\
&\leq P(C_m)P(\hat C_{m,n})+\phi([\frac 12(k-D)]-1)P(C_m).\nonumber
\end{eqnarray}

To make the estimate (\ref{4.4}) suitable for our purposes we recall that according to Lemma 3.1 in \cite{KY}
there exists $\al>0$ such that any cylinder set $C$ defined on an interval of integers $\La=[l,r]$ satisfies
\begin{equation}\label{4.5}
P(C)\leq e^{-\al(r-l)},
\end{equation}
and so
\begin{equation}\label{4.6}
P(\hat C_{m,n})\leq\exp(-\al([\frac 12(k-D)]-1).
\end{equation}
In addition to (\ref{4.4}) we can write also
\begin{equation}\label{4.7}
P(C_m)P(C_n)\leq e^{-\al(r_n-l_n)}P(C_m)\leq e^{-\al(k-D)}P(C_m)
\end{equation}
where we used that by (\ref{4.1}),
\[
r_n-l_n\geq r_n-s_{m,n}+1=r_n+q(n)-r_m-q(m)>k-D.
\]
Observe that by Assumption \ref{ass2.1} there exists at most $K(D+1)$ of $n$'s for which $q(n)-q(m)=k\leq D$,
and so by (\ref{4.1}) the second inequality in (\ref{4.2}) may fail only for at most $K(D+1)$ of $n$'s. For
such $n$'s we use the trivial estimate
\begin{equation}\label{4.8}
|P(T^{-q(m)}C_m\cap T^{-q(n)}C_n)-P(C_m)P(C_n)|\leq P(C_m).
\end{equation}

Now if
\begin{equation}\label{4.9}
l_n+q(n)>r_m+q(m)
\end{equation}
then by the definition of the $\phi$-dependence coefficient we can write by (\ref{4.1}) that
\begin{eqnarray}\label{4.10}
&|P(T^{-q(m)}C_m\cap T^{-q(n)}C_n)-P(C_m)P(C_n)|\\
&\leq\phi(l_n+q(n)-r_m-q(m))P(C_m)\leq\phi(k-D-(r_n-l_n))P(C_m)\nonumber
\end{eqnarray}
but this may not suffice for our purposes when $r_n-l_n$ is large. In this case we proceed as in
(\ref{4.4}), (\ref{4.6}) and (\ref{4.7}) where we take $\hat C_n=[a_{t_n},a_{t_n+1},...,a_{r_n}]$
with $t_n=l_n+[\frac 12(r_n-l_n)]+1$. Then
\[
t_n+q(n)-r_m-q(m)>[\frac 12(r_n-l_n)]+1\quad\mbox{and}\quad r_n-t_n\geq [\frac 12(r_n-l_n)]-1,
\]
and so
\begin{eqnarray}\label{4.11}
&P(T^{-q(m)}C_m\cap T^{-q(n)}C_n)\leq P(T^{-q(m)}C_m\cap T^{-q(n)}\hat C_n)\\
&\leq P(C_m)P(\hat C_n)+\phi([\frac 12(r_n-l_n)]+1)P(C_m)\nonumber\\
&\leq\big( e^{-\al([\frac 12(r_n-l_n)]-1)}+\phi([\frac 12(r_n-l_n)])\big) P(C_m).\nonumber
\end{eqnarray}
Thus, when (\ref{4.9}) holds true we use (\ref{4.10}) if $r_n-l_n\leq\frac {k-D}2$ and (\ref{4.11})
when $r_n-l_n>\frac {k-D}2$. In both cases we will obtain the estimate
\begin{eqnarray}\label{4.12}
&|P(T^{-q(m)}C_m\cap T^{-q(n)}C_n)-P(C_m)P(C_n)|\\
&\leq \big( e^{-\al([\frac 14(k-D)]-1)}+\phi([\frac 14(k-D)])\big) P(C_m).\nonumber
\end{eqnarray}

Finally, taking into account that $q(n)-q(m)=k\leq D$ can occur only for at most $K(D+1)$ of $n$'s
and for each $k$ the equality $q(n)-q(m)=k$ may hold true for at most $K$ of $n$'s we conclude from
(\ref{4.4}), (\ref{4.6})--(\ref{4.8}), (\ref{4.12}) and from the summability of the coefficient $\phi$
that for any $m=M,M+1,...,N$,
\begin{equation}\label{4.13}
\sum_{n=M}^N|P(T^{-q(m)}C_m\cap T^{-q(n)}C_n)-P(C_m)P(C_n)|\leq cP(C_m)
\end{equation}
for some constant $c>0$ independent of $M$ and $N$. Summing in $m$ between $M$ and $N$ we conclude that
the condition (\ref{2.12}) of Theorem \ref{thm2.4} is satisfied with $\Gam_n=T^{-q(n)}C_n$, and so assuming
(\ref{2.10}) we obtain (\ref{2.11}) completing the proof of Theorem \ref{thm2.3}(i).  \qed

\section{Proof of Theorem 2.3($ii$)}\label{sec5}\setcounter{equation}{0}

Observe that if $\del(n,m)=k$, $n>m\geq 0$ and the pair $n,m$ does not belong to the exceptional set
$\cN$ having cardinality at most $K$ then by Assumption \ref{ass2.1}(ii) for some $i_0,j_0\leq\ell$,
\begin{equation}\label{5.1}
q_{j_0}=\max_{1\leq j\leq\ell}q_j(n)\geq q_{i_0}(m)+k=\max_{1\leq i\leq\ell}q_i(m)+k.
\end{equation}
Let $C_m$ and $C_n$ be cylinder sets defined on $\La_m=[l_m,r_m]$ and $\La_n=[l_n,r_n]$, respectively.
Since $C_m$ is $D$-nested in $C_n$, $r_m\leq r_n+D$, and so by (\ref{5.1}),
\begin{equation}\label{5.2}
r_m+q_{i_0}(m)\leq r_n+q_{j_0}(n)-k+D.
\end{equation}

Assume first that
\begin{equation}\label{5.3}
l_n+q_{j_0}(n)\leq r_m+q_{i_0}(m)\quad{and}\quad r_n+q_{j_0}(n)>r_m+q_{i_0}(m).
\end{equation}
Let $C_n=[a_{l_n},a_{l_n+1},...,a_{r_n}]$ and $\hat C_{m,n}=[a_{s_{m,n}},a_{s_{m,n}+1},...,a_{r_n}]$ where
\begin{equation}\label{5.4}
s_{m,n}=l_n+(r_m+q_{i_0}(m)-l_n-q_{j_0})+1=r_m+q_{i_0}(m)-q_{j_0}(n)+1,
\end{equation}
and so $\hat C_{m,n}$ is defined on the interval $[s_{m,n},r_n]$ of the length
\begin{equation}\label{5.5}
r_n-s_{m,n}+1=r_n+q_{j_0}(n)-r_m-q_{i_0}(m)\geq k-D
\end{equation}
where the last inequality follows from (\ref{5.2}). Hence, by the definition of the $\psi$-dependence coefficient
\begin{eqnarray}\label{5.6}
&P\big(\cap_{i=1}^\ell(T^{-q_i(m)}C_m^{(i)}\cap T^{-q_i(n)}C_n^{(i)})\big)\\
&\leq P\big(\cap_{i=1}^\ell(T^{-q_i(m)}C_m^{(i)}\cap T^{-q_{j_0}(n)}\hat C_{m,n}^{(j_0)})\big)\nonumber\\
&\leq (1+\psi(1))P(\cap_{i=1}^\ell T^{-q_i(m)}C_m^{(i)})P( T^{-q_{j_0}(n)}\hat C_{m,n}^{(j_0)})\nonumber\\
&\leq (1+\psi(1))e^{-\al(k-D)}P(\cap_{i=1}^\ell T^{-q_i(m)}C_m^{(i)})\nonumber
\end{eqnarray}
where $\hat C_{m,n}^{(j_0)}$ is constructed as above with $C_n=C_n^{(j_0)}$.

We can write also that
\begin{eqnarray}\label{5.7}
&P(\cap_{i=1}^\ell T^{-q_i(m)}C_m^{(i)})P(\cap_{i=1}^\ell T^{-q_i(n)}C_n^{(i)})\leq P(C_n^{(1)})P(\cap_{i=1}^\ell T^{-q_i(m)}C_m^{(i)})\\
&\leq e^{-\al(r_n-l_n)}P(\cap_{i=1}^\ell T^{-q_i(m)}C_m^{(i)}).\nonumber
\end{eqnarray}
Since $r_n-l_n\geq r_n-s_{m,n}+1\geq k-D$, it follows that under the condition (\ref{5.3}),
\begin{eqnarray}\label{5.8}
&|P\big(\cap_{i=1}^\ell(T^{-q_i(m)}C_m^{(i)}\cap T^{-q_i(n)}C_n^{(i)})\big)\\
&-P(\cap_{i=1}^\ell T^{-q_i(m)}C_m^{(i)})P(\cap_{i=1}^\ell T^{-q_i(n)}C_n^{(i)})|\nonumber\\
&\leq (1+\psi(1))e^{-\al(k-D)}P(\cap_{i=1}^\ell T^{-q_i(m)}C_m^{(i)}).\nonumber
\end{eqnarray}

On the other hand, if
\begin{equation}\label{5.9}
l_n+q_{j_0}(n)>r_m+q_{i_0}(m),
\end{equation}
then by the definition of the $\psi$-dependence coefficient we obtain similarly to the above that
\begin{eqnarray}\label{5.10}
&|P\big(\cap_{i=1}^\ell(T^{-q_i(m)}C_m^{(i)}\cap T^{-q_i(n)}C_n^{(i)})\big)\\
&-P(\cap_{i=1}^\ell T^{-q_i(m)}C_m^{(i)})P(\cap_{i=1}^\ell T^{-q_i(n)}C_n^{(i)})|\nonumber\\
&\leq (1+\psi(l_n+q_{j_0}(n)-r_m-q_{i_0}(m)))P(\cap_{i=1}^\ell T^{-q_i(m)}C_m^{(i)})P(C_n^{(1)})\nonumber\\
&\leq (1+\psi(1))e^{-\al(r_n-l_n)}P(\cap_{i=1}^\ell T^{-q_i(m)}C_m^{(i)}).\nonumber
\end{eqnarray}

Let a number $d_0\geq 1$ be such that
\begin{equation}\label{5.11}
\psi(d_0)<2^{1/\ell}-1\quad\mbox{and}\quad k-(r_n-l_n+2D)>d_0.
\end{equation}
Since $r_n-l_n\geq r_m-l_m-2D$ by $D$-nesting, it follows by (\ref{4.5}) and Lemma 3.3 from \cite{KR}
that
\begin{eqnarray}\label{5.12}
&|P\big(\cap_{i=1}^\ell (T^{-q_i(m)}C_m^{(i)}\cap T^{-q_i(n)}C_n^{(i)})\big)\\
&-P(\cap_{i=1}^\ell T^{-q_i(m)}C_m^{(i)})P(\cap_{i=1}^\ell T^{-q_i(n)}C_n^{(i)})|\nonumber\\
&\leq 2^{2\ell+2}\psi(k-\max(r_n-l_n,r_m-l_m))\nonumber\\
&\times (2-(1+\psi(k-\max(r_n-l_n,r_m-l_m))^\ell)^{-2}P(\cap_{i=1}^\ell T^{-q_i(m)}C_m^{(i)})\nonumber\\
&\times P(\cap_{i=1}^\ell T^{-q_i(n)}C_n^{(i)})\leq 2^{2\ell+2}\psi(k-(r_n-l_n+2D))(2-(1+\psi(d_0))^\ell)^{-2}\nonumber\\
&\times e^{-\al(r_n-l_n)}P(\cap_{i=1}^\ell T^{-q_i(m)}C_m^{(i)}).\nonumber
\end{eqnarray}

Since the cardinality of $\cN$ does not exceed $K$ we have
\begin{eqnarray}\label{5.13}
&\sum_{(n,m)\in\cN}|P\big(\cap_{i=1}^\ell (T^{-q_i(m)}C_m^{(i)}\cap T^{-q_i(n)}C_n^{(i)})\big)\\
&-P(\cap_{i=1}^\ell T^{-q_i(m)}C_m^{(i)})P(\cap_{i=1}^\ell T^{-q_i(n)}C_n^{(i)})|\nonumber\\
&\leq KP(\cap_{i=1}^\ell T^{-q_i(m)}C_m^{(i)}).\nonumber
\end{eqnarray}
Next, we estimate now the remaining sum
\begin{eqnarray}\label{5.14}
&\sum_{n>m,(n,m)\not\in\cN}|P\big(\cap_{i=1}^\ell (T^{-q_i(m)}C_m^{(i)}\cap T^{-q_i(n)}C_n^{(i)})\big)\\
&-P(\cap_{i=1}^\ell T^{-q_i(m)}C_m^{(i)})P(\cap_{i=1}^\ell T^{-q_i(n)}C_n^{(i)})|.\nonumber
\end{eqnarray}
For the part of the sum in $n$'s satisfying (\ref{5.3}) we apply the inequality (\ref{5.8}) which yields the contribution
to the total sum estimated using (\ref{3.4}) by
\begin{eqnarray}\label{5.15}
&2K\ell^2(1+\psi(1))P(\cap_{i=1}^\ell T^{-q_i(m)}C_m^{(i)})\sum_{k=0}^\infty e^{-\al(k-D)}\\
&=2K\ell^2e^{\al D}(1+\psi(1))(1-e^{-\al})^{-1}P(\cap_{i=1}^\ell T^{-q_i(m)}C_m^{(i)}).\nonumber
\end{eqnarray}

For the parts of the sum (\ref{5.14}) which correspond to $n$'s satisfying (\ref{5.9}) but not (\ref{5.11})
we obtain that
\begin{equation}\label{5.17}
e^{-\al(r_n-l_n)}\leq e^{-\al k}e^{-\al(2D-d_0)},
\end{equation}
and so taking into account (\ref{3.4}) the summation in (\ref{5.14}) over $n$'s satisfying (\ref{5.9}) can be
estimated by
\begin{eqnarray}\label{5.18}
&2K\ell^2(1+\psi(1))e^{-\al(2D-d_0)}P(\cap_{i=1}^\ell T^{-q_i(m)}C_m^{(i)})\sum_{k=0}^\infty e^{-\al k}\\
&=2K\ell^2(1+\psi(1))e^{-\al(2D-d_0)}(1-e^{-\al})P(\cap_{i=1}^\ell T^{-q_i(m)}C_m^{(i)})\nonumber.
\end{eqnarray}

It remains to estimate the part of the sum (\ref{5.14}) which corresponds to $n$'s satisfying (\ref{5.11})
where we use (\ref{5.12}). We observe that
\begin{eqnarray}\label{5.19}
&\psi(k-(r_n-l_n+2D))e^{-\al(r_n-l_n)}\\
&=e^{2\al D}\psi(k-(r_n-l_n+2D))e^{-\al(r_n-l_n+2D)}\nonumber\\
&\leq e^{2\al D}\max(\psi([k/2]), \psi(1)e^{-\al[k/2]})\leq e^{2\al D}(\psi([k/2])+\psi(1)e^{-\al[k/2]})\nonumber
\end{eqnarray}
since either $r_n-l_n+2D\geq k/2$ or $k-(r_n-l_n+2D)\geq k/2$. Both summands in the right hand side of (\ref{5.19})
are summable in $k$ (the first one by the assumption) which gives an estimate for the part of the sum (\ref{5.14})
corresponding to $n$'s satisfying (\ref{5.11}) in the form
\begin{equation}\label{5.20}
cP(\cap_{i=1}^\ell T^{-q_i(m)}C_m^{(i)})
\end{equation}
where $c>0$ does not depend on $m$. By estimates (\ref{5.8}), (\ref{5.12}), (\ref{5.13}), (\ref{5.15}) and
(\ref{5.18})--(\ref{5.20}) above we conclude that the whole sum consisting of the part appearing in (\ref{5.13})
plus the part displayed by (\ref{5.14}) can be estimated by the expression (\ref{5.20}) with another constant
$c>0$ independent of $m$. It follows that there exists $\tilde c>0$ such that for all $N>M\geq 1$,
\begin{eqnarray}\label{5.21}
&\sum_{n,m=M}^N|P\big(\cap_{i=1}^\ell (T^{-q_i(m)}C_m^{(i)}\cap T^{-q_i(n)}C_n^{(i)})\big)\\
&-P(\cap_{i=1}^\ell T^{-q_i(m)}C_m^{(i)})P(\cap_{i=1}^\ell T^{-q_i(n)}C_n^{(i)})|\nonumber\\
&\leq\tilde c\sum_{m=M}^NP(\cap_{i=1}^\ell T^{-q_i(m)}C_m^{(i)}).\nonumber
\end{eqnarray}
If, in addition,
\begin{equation}\label{5.22}
\sum_{m=1}^\infty P(\cap_{i=1}^\ell T^{-q_i(m)}C_m^{(i)})=\infty
\end{equation}
then by Theorem \ref{thm2.4} we obtain that with probability one
\begin{equation}\label{5.23}
\frac {\sum_{n=1}^N(\prod_{i=1}^\ell\bbI_{C_n^{(i)}}\circ T^{q_i(n)})}{\sum_{n=1}^NP(\cap_{i=1}^\ell T^{-q_i(n)}C_n^{(i)})}\to 1\,\,\mbox{as}\,\, N\to\infty.
\end{equation}

It remains to show that under the condition (\ref{2.10}) with probability one,
\begin{equation}\label{5.24}
\frac {\sum_{n=1}^NP(\cap_{i=1}^\ell T^{-q_i(n)}C_n^{(i)})}{\sum_{n=1}^N\prod_{i=1}^\ell P(C_n^{(i)})}\to 1\quad\mbox{as}\quad N\to\infty.
\end{equation}
Observe again that
\begin{equation}\label{5.25}
P(\cap_{i=1}^\ell T^{-q_i(n)}C_n^{(i)})\leq P(C_n^{(1)})\leq e^{-\al(r_n-l_n)}.
\end{equation}
Next, we split the sum in the left hand side of (\ref{5.22}) into two sums
\begin{eqnarray*}
&S_1=\sum_{n:\,(r_n-l_n)\leq\frac 2\al\ln n}P(\cap_{i=1}^\ell(T^{-q_i(n)}C_n^{(i)})\\
&\mbox{and}\,\,\,S_2=\sum_{n:\,(r_n-l_n)>\frac 2\al\ln n}P(\cap_{i=1}^\ell(T^{-q_i(n)}C_n^{(i)}).
\end{eqnarray*}
By (\ref{5.25}),
\[
S_2\leq\sum_{n=1}^\infty n^{-2}<\infty\,\,\mbox{and also}\,\,\sum_{n:\,(r_n-l_n)>\frac 2\al\ln n}\prod_{i=1}^\ell P(C_n^{(i)})<\infty.
\]
Hence, it suffices to show that under the condition (\ref{2.10})  with probability one,
\begin{equation}\label{5.26}
\frac {\sum_{n\leq N:\,(r_n-l_n)\leq\frac 2\al\ln n}P(\cap_{i=1}^\ell T^{-q_i(n)}C_n^{(i)})}{\sum_{n\leq N:\,n:\,(r_n-l_n)\leq\frac 2\al\ln n}
\prod_{i=1}^\ell P(C_n^{(i)})}\to 1\quad\mbox{as}\quad N\to\infty.
\end{equation}

Set $q(n)=\min_{i\ne j}|q_i(n)-q_j(n)|$. Observe that by Assumption \ref{ass2.1}(i) for each $k$,
\begin{equation}\label{5.27}
\#\{ n:\, q(n)=k\}\leq K\ell^2.
\end{equation}
Consider first $n$'s satisfying
\begin{equation}\label{5.28}
q(n)\leq r_n-l_n.
\end{equation}
In this case by (\ref{5.25}),
\begin{equation}\label{5.29}
P(\cap^\ell_{i=1}T^{-q_i(n)}C_n^{(i)})\leq e^{-\al q(n)}
\end{equation}
and relying on (\ref{5.27}) we conclude that
\[
\sum_{n:\, q(n)\leq r_n-l_n}P(\cap^\ell_{i=1}T^{-q_i(n)}C_n^{(i)})\leq K\ell^2\sum_{k=0}^\infty e^{-\al k}=K\ell^2(1-e^{-\al})^{-1}
\]
and the same estimate holds true for $\sum_{n:\, q(n)\leq r_n-l_n}\prod_{i=1}^\ell P(C_n^{(i)})$. Hence, the sum over such $n$'s does not influence
the asymptotical behavior in (\ref{5.24}) and (\ref{5.26}) since the denominators there tend to $\infty$.

It remains to consider the sums over $n$'s satisfying
\begin{equation}\label{5.30}
q(n)>r_n-l_n.
\end{equation}
In this case we can apply Lemma 3.2 from \cite{KR} to obtain that
\begin{eqnarray}\label{5.31}
&|P(\cap^\ell_{i=1}T^{-q_i(n)}C_n^{(i)})-\prod_{i=1}^\ell P(C_n^{(i)})|\\
&\leq\big((1+\psi(q(n)-(r_n-l_n))^\ell-1\big)\prod_{i=1}^\ell P(C_n^{(i)})\nonumber\\
&\leq \big((1+\psi(q(n)-(r_n-l_n))^\ell-1\big)e^{-\ell\al(r_n-l_n)}.\nonumber
\end{eqnarray}

Now observe that either $r_n-l_n$ or $q(n)-(r_n-l_n)$ is greater or equal to $\frac 12q(n)$. Denote by $\cN_1$ the set of $n$'s
for which $r_n-l_n\geq \frac 12q(n)$ and by $\cN_2$ the set of $n$'s for which $q(n)-(r_n-l_n)\geq\frac 12q(n)$. Taking into account
(\ref{5.27}) and (\ref{5.30}) we obtain that
\begin{eqnarray}\label{5.32}
&\sum_{n\in\cN_1}\big((1+\psi(q(n)-(r_n-l_n))^\ell-1\big)e^{-\ell\al(r_n-l_n)}\\
&\leq((1+\psi(1))^\ell-1)\sum_{n\in\cN_1}e^{-\frac 12\ell\al q(n)}\nonumber\\
&\leq K\ell^2((1+\psi(1))^\ell-1)\sum_{k=0}^\infty e^{-\frac 12\ell\al k}\nonumber\\
&=K\ell^2((1+\psi(1))^\ell-1)(1-e^{-\frac 12\ell\al})^{-1}<\infty.\nonumber
\end{eqnarray}
Next, taking into account that $\psi(k)$ is summable we see that
\begin{eqnarray}\label{5.33}
&\sum_{n\in\cN_2}\big((1+\psi(q(n)-(r_n-l_n))^\ell-1\big)e^{-\ell\al(r_n-l_n)}\\
&\leq\sum_{n\in\cN_2}\big((1+\psi(\max(1,[\frac 12q(n)]))^\ell-1\big)\nonumber\\
&\leq 2K\ell^2\sum_{k=1}^\infty((1+\psi(k))^\ell-1)=2K\ell^2\sum_{k=1}^\infty\sum_{m=1}^\ell{\ell\choose m}(\psi(k))^m<\infty.\nonumber
\end{eqnarray}
Hence,
\begin{equation}\label{5.34}
|\sum_{n=1}^\infty\big( P(\cap^\ell_{i=1}T^{-q_i(n)}C_n^{(i)})-\prod_{i=1}^\ell P(C_n^{(i)})\big)|<\infty
\end{equation}
and since $\sum_{n=1}^\infty\prod_{i=1}^\ell P(C_n^{(i)})=\infty$, we obtain (\ref{5.26}), and so (\ref{5.24}), as well, completing the proof of
Theorem \ref{thm2.3}(ii).   \qed

\section{Asymptotics of maximums of logarithmic distance functions}\label{sec6}\setcounter{equation}{0}

In this section we will prove Theorem \ref{thm2.5}. Let $\tilde\om^{(j)}=(\tilde\om_i^{(j)})_{i\in\bbZ}\in\Om$ and $C_n(\tilde\om^{(j)}),\, j=1,...,\ell,\, n=1,2,...$ be a sequence of cylinder sets such that
\[
C_n(\tilde\om^{(j)})=\{\om=(\om_i)_{i\in\bbZ}\in\Om:\,\om_i=\tilde\om_i^{(j)}\,\,\mbox{ provided}\,\, |i|\leq r_n\}
\]
where $r_n\uparrow\infty$ as $n\uparrow\infty$ is a sequence of integers. Observe that by the Shannon--McMillan--Breiman
theorem (see, for instance, \cite{Pe}) for almost all $\tilde\om\in\Om$,
\begin{equation}\label{6.1}
\lim_{n\to\infty}\frac 1{2r_n}\ln P(C_n(\tilde\om))=-h
\end{equation}
where $h$ is the Kolmogorov--Sinai entropy of the shift $T$ with respect to $P$ since the latter measure is ergodic whether
we assume $\phi$ or $\psi$-mixing.

Now suppose that
\begin{equation}\label{6.2}
\sum_{n=1}^\infty \prod_{i=1}^\ell P(C_n(\tilde\om^{(i)}))<\infty.
\end{equation}
It follows from (\ref{5.33}) that (\ref{6.2}) implies also
\begin{equation}\label{6.3}
\sum_{n=1}^\infty P(\cap_{i=1}^\ell T^{-q_i(n)}C_n(\tilde\om^{(i)}))<\infty
\end{equation}
which is, of course, a tautology if $\ell=1$. It follows from the first Borel--Cantelli lemma that for almost all $\om\in\Om$
only finitely many events $\{ T^{q_i(n)}\om\in C_n(\tilde\om^{(i)}),\, i=1,...,\ell\}$ can occur. But if the latter event does not hold true then
\[
T^{q_j(n)}\not\in C_n(\tilde\om^{(j)})\,\,\,\mbox{for some}\,\,\, 1\leq j\leq\ell,
\]
and so
\begin{equation}\label{6.4}
d(T^{q_j(n)}\om,\tilde\om^{(j)})>e^{-\gam r_n}\,\,\mbox{i.e.}\,\,\Phi_{\tilde\om^{(j)}}(T^{q_j(n)}\om)<\gam r_n
\end{equation}
where the distance $d(\cdot,\cdot)$ and the function $\Phi$ were defined in (\ref{2.14}) and (\ref{2.15}).
It follows that in this case there exists $N_{\bm{\tilde\om}},\,\bm{\tilde\om}=(\tilde\om^{(1)},...,\tilde\om^{(\ell)})$ finite with probability one and such that for
all $N>N_{\bm{\tilde\om}}(\om)$,
\[
M_{N,\bm{\tilde\om}}(\om)<\gam r_N,
\]
where $M_{N,\bm{\tilde\om}}(\om)$ was defined in (\ref{2.15}). Hence,
\begin{equation}\label{6.5}
\limsup_{N\to\infty}\frac {M_{N,\bm{\tilde\om}}}{\ln N}\leq\gam\limsup_{N\to\infty}\frac {r_N}{\ln N}\,\,\mbox{a.s.}
\end{equation}

Next, assume that
\begin{equation}\label{6.6}
\cE_{N,\tilde\om}=\sum_{n=1}^N\prod_{i=1}^\ell P(C_n(\tilde\om^{(i)})\to \infty\,\,\mbox{as}\,\, N\to\infty
\end{equation}
which by (\ref{5.33}) implies also that
\begin{equation}\label{6.7}
\sum_{n=1}^\infty P(\cap_{i=1}^\ell T^{-q_i(n)}C_n(\tilde\om^{(i)}))=\infty.
\end{equation}
Set
\[
L_{n,\bm{\tilde\om}}(\om)=\max\{ m\leq n:\, T^{q_i(m)}\om\in C_m(\tilde\om^{(i)})\,\,\mbox{for}\,\, i=1,...,\ell\}.
\]
It follows from Theorem \ref{thm2.3} that under (\ref{6.6}) for almost all $\om\in\Om$,
\[
L_{n,\bm{\tilde\om}}(\om)\to\infty\,\,\,\mbox{as}\,\,\, n\to\infty.
\]
Observe also that
\begin{equation}\label{6.8}
S_N(\om)=\sum_{n=1}^N(\prod_{i=1}^\ell\bbI_{C_n(\tilde\om^{(i)})}\circ T^{q_i(n)}(\om))=S_{L_{n,\bm{\tilde\om}}(\om)}.
\end{equation}
By (\ref{4.13}), (\ref{5.20}) and (\ref{5.33}) we can use (\ref{2.13}) which yields that for almost all $\om\in\Om$,
\begin{equation}\label{6.9}
0\leq\cE_{N,\bm{\tilde\om}}-\cE_{L_{n,\bm{\tilde\om}},\bm{\tilde\om}}\leq O(\cE_{N,\bm{\tilde\om}}^{1/2}\ln^{\frac 32+\ve}\cE_{N,\bm{\tilde\om}}),
\end{equation}
and so for almost all $\om$,
\begin{equation}\label{6.10}
\lim_{N\to\infty}\frac {\cE_{L_{n,\bm{\tilde\om}}(\om),\bm{\tilde\om}}}{\cE_{N,\bm{\tilde\om}}}=1.
\end{equation}

Next, observe that if $m=L_{n,\bm{\tilde\om}}(\om)$ then for each $i=1,...,\ell$,
\[
d(T^{q_i(m)}\om,\tilde\om^{(i)})\leq e^{-\gam r_m}\,\,\mbox{i.e.}\,\,\Phi_{\tilde\om^{(i)}}(T^{q_i(m)}\om)\geq\gam r_m,
\]
and so $M_{m,\bm{\tilde\om}}(\om)\geq\gam r_m$. It follows that
\begin{equation}\label{6.11}
M_{N,\bm{\tilde\om}}(\om)\geq M_{L_{N,\bm{\tilde\om}}(\om)}(\om)\geq\gam r_{L_{N,\bm{\tilde\om}}(\om)},
\end{equation}
and so
\begin{equation}\label{6.12}
\liminf_{N\to\infty}\frac {M_{N,\bm{\tilde\om}}(\om)}{\ln N}\geq\gam(\liminf_{N\to\infty}\frac {r_N}{\ln N})\liminf_{N\to\infty}
\frac {\ln L_{N,\bm{\tilde\om}}(\om)}{\ln N}.
\end{equation}

Next, in order to complete the proof of Theorem \ref{thm2.5}, we will choose sequences $r_n,\, n=1,2,...$ for appropriate upper
and lower bounds. For the upper bound we will take $r_n=[\frac {1+\del}{2\ell h}\ln n]$ for some $\del>0$. Then by (\ref{6.1}) for almost all $\tilde\om^{(1)},...,\om^{(\ell)}\in\Om$,
\[
\ln\prod_{i=1}^\ell P(C_n(\tilde\om^{(i)}))\sim -(1+\del)\ln n\quad\mbox{as}\quad n\to\infty,
\]
and so the series (\ref{6.2}) converges as needed. Substituting such $r_N$'s to (\ref{6.5}) and letting $\del\to 0$ we obtain
\begin{equation}\label{6.13}
\limsup_{N\to\infty}\frac {M_{N,\bm{\tilde\om}}}{\ln N}\leq\frac \gam{2\ell h}\quad\mbox{a.s.}
\end{equation}

Now we deal with the lower bound choosing $r_n=[\frac {1-\del}{2\ell h}\ln n]$. Then by (\ref{6.1}) for almost all
 $\tilde\om^{(1)},...,\om^{(\ell)}\in\Om$ as $n\to\infty$,
\begin{equation}\label{6.14}
\ln\prod_{i=1}^\ell P(C_n(\tilde\om^{(i)}))\sim -(1-\del)\ln n,
\end{equation}
and so the series (\ref{6.6}) diverges as needed. For such $r_N$'s we have that
\begin{equation}\label{6.15}
\liminf_{N\to\infty}\frac {r_N}{\ln N}=\frac {1-\del}{2\ell h}
\end{equation}
and letting $\del\to 0$ the proof of Theorem \ref{thm2.5} will be completed by (\ref{6.12}), (\ref{6.13}) and
(\ref{6.15}) once we show that for
almost all $\om\in\Om$,
\begin{equation}\label{6.16}
\liminf_{N\to\infty}\frac {\ln L_{N,\bm{\tilde\om}}(\om)}{\ln N}=1.
\end{equation}

By (\ref{6.14}) there exists a random variable $n(\bm{\tilde\om})<\infty$ a.s. such that if $n\geq n(\bm{\tilde\om})$ then
\begin{equation}\label{6.17}
n^{-(1-\frac 34\del)}\leq\prod_{i=1}^\ell P(C_n(\tilde\om^{(i)})) \leq n^{-(1-\frac 43\del)}.
\end{equation}
If $L_{N,\tilde\om}(\om)\geq n(\om)$ then we obtain from (\ref{6.9}) and (\ref{6.17}) that
\begin{eqnarray}\label{6.18}
&\frac 4{3\del}(N^{\frac 34\del}-(L_{N,\bm{\tilde\om}}(\om)+1)^{\frac 34\del})\leq\sum_{n=L_{N,\bm{\tilde\om}}(\om)+1}^Nn^{-(1-\frac 34\del)}\\
&\leq O\big( (n(\om)+\sum_{n=n(\om)}^Nn^{-(1-\frac 43\del)})^{1/2}\ln^{\frac 32+\ve}(n(\om)+\sum_{n=n(\om)}^Nn^{-(1-\frac 43\del)})\big)\nonumber\\
&\leq O\big(n(\om)+\frac 3{4\del}N^{\frac 43\del})^{1/2}\ln^{\frac 32+\ve}(n(\om)+\frac 3{4\del}N^{\frac 43\del})\big).\nonumber
\end{eqnarray}
Dividing these inequalities by $N^{\frac 34\del}$, letting $N\to\infty$ and taking into account that $N\geq L_{N,\bm{\tilde\om}}(\om)$ by the definition, we see that
\[
\frac {L_{N,\bm{\tilde\om}}(\om)}N\to 1,\,\,\mbox{and so}\,\, \ln N-\ln L_{N,\tilde\om}(\om)\to 0\,\,\mbox{a.s. as}
\,\, N\to\infty
\]
implying (\ref{6.16}) and completing the proof of Theorem \ref{thm2.5}.   \qed

\section{Asymptotics of hitting times}\label{sec7}\setcounter{equation}{0}

In this section we will prove Theorem \ref{thm2.6} deriving first that for $P\times P$-almost all pairs
$(\om,\tilde\om)$,
\begin{equation}\label{7.1}
\liminf_{n\to\infty}\frac 1n\ln\tau_{C_n(\tilde\om)}\geq 2\ell h.
\end{equation}
Let $\la k\leq n\leq\la(k+1)$ for some $\la>0$. Then
\[
\frac {\ln\tau_{C_{\la k}(\tilde\om)}}{\la(k+1)}\leq\frac {\ln\tau_{C_n(\tilde\om)}}{n}\leq
\frac {\ln\tau_{C_{\la( k+1)}(\tilde\om)}}{\la k},
\]
and so
\begin{equation}\label{7.2}
\liminf_{n\to\infty}\frac {\ln\tau_{C_n(\tilde\om)}}{n}=\liminf_{k\to\infty}
\frac {\ln\tau_{C_{\la k}(\tilde\om)}}{\la k}
\end{equation}
where we alert the reader that the definition of the cylinder $C_n\tilde\om)$ here agrees with
the corresponding definition in Section \ref{sec6} provided $r_n=n$ there.

Next, assume that $\la>(2\ell h)^{-1}$ and set
\[
I_k(\tilde\om)=\cup_{j=1}^{e^k}\cap_{i=1}^\ell T^{-q_i(j)}C_{\la k}(\tilde\om).
\]
Then
\begin{equation}\label{7.3}
P(I_k(\tilde\om))\leq\sum_{j=1}^{e^k}P(\cap_{i=1}^\ell T^{-q_i(j)}C_{\la k}(\tilde\om))
\end{equation}
and we are going to show that for $P$-almost all $\tilde\om$,
\begin{equation}\label{7.4}
\sum_{k=1}^\infty\sum_{j=1}^{e^k}P(\cap_{i=1}^\ell T^{-q_i(j)}C_{\la k}(\tilde\om))<\infty.
\end{equation}
Indeed, applying the Shannon-McMillan-Breiman theorem we obtain that for $P$-almost all $\tilde\om$
and each $\ve>0$ there exists $k(\ve,\tilde\om)$ such that if $k\geq k(\ve,\tilde\om)$ then
\begin{equation}\label{7.5}
P(C_{\la k}(\tilde\om))\leq\exp(-k(2\la h-\ve)).
\end{equation}
When $\ell=1$ we employ (\ref{7.5}) for $k\geq k(\ve,\tilde\om)$ and (\ref{4.5}) for
$k< k(\ve,\tilde\om)$ which yields the estimate of the left hand side of (\ref{7.4}) by
\begin{equation}\label{7.6}
\sum_{1\leq k\leq k(\ve,\tilde\om)}e^ke^{-\al(2\la k-1)}+\sum_{k=1}^\infty e^{-k(2\la h-\ve-1)}.
\end{equation}
The first sum in (\ref{7.6}) contains finitely many terms, and so it is bounded, while the
second sum in (\ref{7.6}) is also bounded since $2\la h-\ve>1$ by the choice of $\la$ provided
$\ve>0$ is small enough.

Next, we will deal with the case $\ell>1$. First, recall the notation
$q(n)=\min_{i\ne j}|q_i(n)-q_j(n)|$ and observe that by (\ref{5.6}),
\begin{equation}\label{7.7}
\#\{ n:\, q(n)\leq 2\la k+2\}\leq 2(\la k+1)K\ell^2.
\end{equation}
Now we split the sum in the left hand side of (\ref{7.4}) into two sums
\begin{eqnarray}\label{7.8}
&S_1=\sum^\infty_{k=1}\sum_{j:j\leq e^k,\, q(j)\leq 2\la k+2}P(\cap_{i=1}^\ell T^{-q_i(j)}
C_{\la k}(\tilde\om))\\
&\leq 2K\ell^2\sum_{k=1}^\infty(\la k+1)P(C_{\la k}(\tilde\om))\leq 2K\ell^2\sum_{k=1}^\infty(\la k+1)
e^{-2\al(\la k-1)}<\infty,\nonumber
\end{eqnarray}
where we we use (\ref{4.5}), and
\begin{eqnarray}\label{7.9}
&S_2=\sum^\infty_{k=1}\sum_{j:j\leq e^k,\, q(j)> 2\la k+2}P(\cap_{i=1}^\ell T^{-q_i(j)}
C_{\la k}(\tilde\om))\\
&\leq k(\ve,\tilde\om)e^{k(\ve,\tilde\om)}+\sum^\infty_{k=k(\ve,\tilde\om)}
\sum_{j:j\leq e^k,\, q(j)> 2\la k+2}P(\cap_{i=1}^\ell T^{-q_i(j)}C_{\la k}(\tilde\om)).\nonumber
\end{eqnarray}
If $q(j)>2\la k+2$ and $k\geq k(\ve,\tilde\om)$ then employing Lemma 3.2 from \cite{KR} and
(\ref{7.5}) above we obtain
\begin{equation}\label{7.10}
P(\cap_{i=1}^\ell T^{-q_i(j)}C_{\la k}(\tilde\om))\leq (1+\psi(1))^\ell(P(C_{\la k}(\tilde\om)))^\ell\leq(1+\psi(1))^\ell\exp(-k(2\la h\ell-\ve\ell))
\end{equation}
where $\psi$ is the dependence coefficient from (\ref{2.3}). For $\ve>0$ small enough
$2\la h\ell-\ve\ell>1$ by the choice of $\la$, and so by (\ref{7.9}) and (\ref{7.10}),
\[
S_2\leq k(\ve,\tilde\om)e^{k(\ve,\tilde\om)}+\sum^\infty_{k=1}\exp(-k(2\la h\ell-\ve\ell-1))
<\infty\]
which together with (\ref{7.8}) yields (\ref{7.4}).

Hence, by the (first) Borel--Cantelli lemma there exists $K(\om)=K(\om,\tilde\om)<\infty$ a.s. 
such that for all $k\geq K(\om)$ there are no events
\[
T^{q_i(j)}\om\in C_{\al k}(\tilde\om)\,\,\,\mbox{for all}\,\,\, i=1,...,\ell\,\,\,\mbox{and some}
\,\,\, 1\leq j\leq e^k.
\]
It follows that for $P$-almost all $\om$ and $k\geq K(\om)$.
\[
\tau_{C_{\al k}(\tilde\om)}(\om)>e^k.
\]
This together with (\ref{7.2}) yields that for $P\times P$-almost all pairs $(\om,\tilde\om)$,
\begin{equation}\label{7.11}
\liminf_{n\to\infty}\frac {\ln\tau_{C_{n}(\tilde\om)}(\om)}n\geq\la^{-1}.
\end{equation}
Since $\la$ can be chosen arbitrarily close to $(2\ell h)^{-1}$ we obtain (\ref{7.1}).

Next, we will prove that for $P\times P$-almost all pairs $(\om,\tilde\om)$,
\begin{equation}\label{7.12}
\limsup_{n\to\infty}\frac {\ln\tau_{C_n(\tilde\om)}(\om)}n\leq 2\ell h.
\end{equation}
Choose $\ve'>\ve>0$ small and $\be>0$ close to $(2\ell h)^{-1}$ so that
\begin{equation}\label{7.13}
\be(2\ell h+\ve)<1\,\,\mbox{and}\,\, \be(2\ell h+\ve')-\frac {1-\be(2\ell h-\ve)}{1-\be(2\ell h+\ve)}>0
\end{equation}
which implies, in particular, that $\be(2\ell h+\ve')>1$.

Set
\[
\Gam=\{(\om,\tilde\om)\in\Om:\,\limsup_{n\to\infty}\frac {\ln\tau_{C_n(\tilde\om)}(\om)}n>2\ell h+\ve'\}.
\]
If $(\om,\tilde\om)\in\Gam$ then for infinitely many $n$'s,
\begin{equation}\label{7.14}
\tau_{C_{\be\ln n}(\tilde\om)}(\om)>n^{\be(2\ell h+\ve')}.
\end{equation}
For $n$'s satisfying (\ref{7.14}),
\[
\om\not\in\cup_{1\leq j\leq n^{\be(2\ell h+\ve')}}\cap_{1\leq i\leq\ell}T^{-q_i(j)}C_{\be\ln n}
\supset\cup_{n\leq j\leq n^{\be(2\ell h+\ve')}}\cap_{1\leq i\leq\ell}T^{-q_i(j)}C_{\be\ln j}
\]
which implies that there exists a sequence $n_k\to\infty$ as $k\to\infty$ such that
\begin{equation}\label{7.15}
\sum_{1\leq j\leq n_k}\prod_{1\leq i\leq\ell}\bbI_{C_{\be\ln j}(\tilde\om)}\circ T^{q_i(j)}(\om)=
\sum_{1\leq j\leq n_k^{\be(2\ell h+\ve')}}\prod_{1\leq i\leq\ell}\bbI_{C_{\be\ln j}(\tilde\om)}\circ T^{q_i(j)}(\om)
\end{equation}
for each $k$.

By the Shannon--McMillan--Breiman theorem there are $\tilde\Om\subset\Om$ with $P(\tilde\Om)=1$ and
a random variable $J$ finite on $\tilde\Om$ such that for any $j\geq J(\tilde\om)$,
\[
j^{-\be(2h+\frac \ve\ell)}=e^{-(2h+\frac \ve\ell)\be\ln j}\leq P(C_{\be\ln j}(\tilde\om))<
e^{-(2h-\frac \ve\ell)\be\ln j}=j^{-\be(2h-\frac \ve\ell)}.
\]
Hence, there are random variables $k_1$ and $k_2$ such that for all $n$ large enough,
\begin{equation}\label{7.16}
k_1(\tilde\om)n^{1-\be(2\ell h+\ve)}\leq\sum_{j=1}^n\big(P(C_{\be\ln j}(\tilde\om))\big)^\ell\leq k_2(\tilde\om)n^{1-\be(2\ell h-\ve)}.
\end{equation}
It follows that for $(\om,\tilde\om)\in\Gam,\,\tilde\om\in\tilde\Om$ and all $k$ large enough
\begin{eqnarray}\label{7.17}
&\frac {\sum_{1\leq j\leq n_k}\big(P(C_{\be\ln j}(\tilde\om))\big)^\ell}
 {\sum_{1\leq j\leq n_k^{\be(2\ell h+\ve')}}\big(P(C_{\be\ln j}(\tilde\om))\big)^\ell}\\
&\leq\frac {k_2(\tilde\om)}{k_1(\tilde\om)}n_k^{(1-\be(2\ell h-\ve)-\be(2\ell h+\ve')(1-\be(2\ell h+\ve))}\to 0\,\,\mbox{as}\,\, k\to\infty
\nonumber\end{eqnarray}
since by the choice of $\ve,\ve'$ and $\be$,
\begin{eqnarray*}
&(1-\be(2\ell h-\ve)-\be(2\ell h+\ve')(1-\be(2\ell h+\ve))\\
&=(1-\be(2\ell h+\ve)(\frac {1-\be(2\ell h-\ve)}{1-\be(2\ell h+\ve)}-\be(2\ell h+\ve')<0.
\end{eqnarray*}
By (\ref{7.15}) we obtain from (\ref{7.17}) that
\begin{eqnarray}\label{7.18}
&\frac {\sum_{1\leq j\leq n_k}\prod_{1\leq i\leq\ell}\bbI_{C_{\be\ln j}(\tilde\om)}
\circ T^{q_i(j)}(\om)}{\sum_{1\leq j\leq n_k}\big(P(C_{\be\ln j}(\tilde\om))\big)^\ell}\\
&\times\frac {\sum_{1\leq j\leq n_k^{\be(2\ell h+\ve')}}\big(P(C_{\be\ln j}(\tilde\om))\big)^\ell}
{\sum_{1\leq j\leq n_k^{\be(2\ell h+\ve')}}\prod_{1\leq i\leq\ell}\bbI_{C_{\be\ln j}(\tilde\om)}\circ T^{q_i(j)}(\om)}\to\infty\,\,\,\mbox{as}\,\,\, k\to\infty.\nonumber
\end{eqnarray}

By (\ref{7.16}) for all $\tilde\om\in\tilde\Om$,
\[
\sum_{1\leq j\leq n}\big(P(C_{\be\ln j}(\tilde\om))\big)^\ell\to\infty\,\,\,\mbox{as}\,\,\, n\to\infty,
\]
and so by Theorem \ref{thm2.3} for $P$-almost all $\om$,
\[
\frac {\sum_{1\leq j\leq n}\prod_{1\leq i\leq\ell}\bbI_{C_{\be\ln j}(\tilde\om)}
\circ T^{q_i(j)}(\om)}{\sum_{1\leq j\leq n}\big(P(C_{\be\ln j}(\tilde\om))\big)^\ell}
\to 1\,\,\,\mbox{as}\,\,\, n\to\infty.
\]
Thus, (\ref{7.18}) can hold true only for a set of pairs $(\om,\tilde\om)$ having $P\times P$-measure zero, and so $P\times P(\Gam)=0$. Since $\ve$ and $\ve'$ can be chosen arbitrarily close to zero,
(\ref{7.12}) follows for $P\times P$-almost all $(\om,\tilde\om)$, which together with (\ref{7.1}) completes the proof of Theorem \ref{thm2.6}.  \qed


\end{document}